\documentclass[11pt,reqno]{amsart}
\usepackage{amsmath,amsthm,amsfonts,amssymb,amscd,amstext}
\usepackage{graphicx}
\usepackage{euscript}
\usepackage{a4wide}
\usepackage{mathpazo} 

\usepackage[shortlabels]{enumitem} 
\numberwithin{equation}{section}

\usepackage[colorlinks=true,backref]{hyperref}

\newcommand{\qand}{\quad\text{and}\quad}

\theoremstyle{plain}
\newtheorem{maintheorem}{Theorem}
\newtheorem{maincorollary}[maintheorem]{Corollary}

\newtheorem{theorem}{Theorem}[section]

\newtheorem{lemma}[theorem]{Lemma}

\theoremstyle{definition}
\newtheorem{remark}[theorem]{Remark}

\newtheorem{conjecture}{Conjecture}

\newcommand{\RR}{{\mathbb R}}

\newcommand{\BB}{{\mathbb B}}

\newcommand{\TT}{{\mathbb T}}
\newcommand{\DD}{{\mathbb D}}
\newcommand{\sS}{{\mathbb S}}

\newcommand{\vfi}{\varphi}

\newcommand{\la}{\lambda}

\renewcommand{\epsilon}{\varepsilon}
\newcommand{\dist}{\operatorname{dist}}

\newcommand{\interior}{\operatorname{int}}

\newcommand{\sing}{\mathrm{Sing}}

\newcommand{\cL}{\EuScript{L}}

\newcommand{\cP}{\EuScript{P}}

\newcommand{\U}{\EuScript{U}}

\newcommand{\cS}{\EuScript{S}}

\newcommand \cR {{\mathcal R}}

\newcommand{\m}{{\rm Leb}\:}

\begin{document}

\title[Physical measures with mixed central behavior] { A
  characterization of physical measures for systems with mixed central
  behavior}

\date{\today}


\thanks{V.A. was partially supported by CNPq-Brazil (grant
  304047/2023-6). L.S. was partially supported by FAPERJ-Funda\c c\~ao
  Carlos Chagas Filho de Amparo \`a Pesquisa do Estado do Rio de
  Janeiro Projects APQ1-E-26/211.690/2021 SEI-260003/015270/2021 and
  JCNE-E-26/200.271/2023 SEI-260003/000640/2023, by Coordena\c c\~ao
  de Aperfei\c coamento de Pessoal de N\'ivel Superior CAPES — Finance
  Code 001 and PROEXT-PG project Dynamic Women - Din\^amicas,
  CNPq-Brazil (grant Projeto Universal 404943/2023-3).}




\author{Vitor Araujo and Luciana Salgado}

\address[V.A.]{Universidade Federal da Bahia,
Instituto de Matem\'atica\\
Av. Adhemar de Barros, S/N , Ondina,
40170-110 - Salvador-BA-Brazil}
\email{vitor.d.araujo@ufba.br \text{or}
  vitor.araujo.im.ufba@gmail.com}
 \urladdr{https://sites.google.com/view/vitor-araujo-ime-ufba}

\address[L.S.]{Universidade Federal do Rio de Janeiro, Instituto de
   Matem\'atica\\
   Avenida Athos da Silveira Ramos 149 Cidade Universit\'aria, P.O. Box 68530, 
   21941-909 Rio de Janeiro-RJ-Brazil }
 \email{lsalgado@im.ufrj.br, lucianasalgado@ufrj.br}
 \urladdr{http://www.im.ufrj.br/~lsalgado}

\begin{abstract}
  We show that the existence of physical measures for $C^\infty$
  smooth instances of certain partially hyperbolic dynamics, both
  continuous and discrete, exhibiting mixed behavior (positive and
  negative Lyapunov exponents) along the central non-uniformly
  hyperbolic multidimensional invariant direction, is equivalent to
  the existence of certain types of ``regular points'' on positive
  volume subsets, including Lyapunov regular points. This encompasses the
  $C^3$ robust class of multidimensional non-hyperbolic attractors
  obtained by Viana, and the $C^1$ robust classes of $3$-sectionally
  hyperbolic wild strange attractors presented by Shilnikov and
  Turaev, providing necessary and sufficient conditions for the
  existence of ergodic hyperbolic physical measures on these and other
  dynamical systems.
\end{abstract}


\maketitle
\tableofcontents

\section{Introduction and statement of results}
\label{sec:intro}

The statistical viewpoint on Dynamical Systems is one of the
cornerstones of most recent developments in dynamics. Given
a map $f$  from a manifold $M$ into itself, a central
concept is that of \emph{physical measure}, a
$f$-invariant probability measure $\mu$ whose
\emph{ergodic basin}
\begin{align*}
  B(\mu)=\Big\{
  x\in M: \lim_{n\to+\infty}\frac1n S_n\varphi(x)=
  \int\varphi\, d\mu\text,\quad \forall \varphi\in C(M,\RR)\Big\}
\end{align*}
having positive \emph{volume} or \emph{Lebesgue measure}, where we write
$S_n \vfi(x) = S_n^f \vfi(x) = \sum_{j=0}^{n-1}
\vfi( f^j x )$ for the standard ergodic sum.
In the continuous time case of a flow $X_t$ generated by the vector
field $X$ on the manifold $M$,
this notion naturally translates to
\begin{align*}
B(\mu)=\Big\{x\in M:
\lim_{T\to+\infty}\frac1T\int_0^T\varphi(X_t(x))=\int\varphi\,
  d\mu, \quad \forall \varphi\in C(M,\RR)\Big\}
\end{align*}
has positive volume.
This kind of measures provides asymptotic information on a
set of trajectories that one expects to be large enough to be
observable in real-world models.

In particular, the existence of these measures implies that the subset
\begin{align*}
  \BB:=\left\{
  x\in M : \Big(\frac1n S_n\delta_x\Big)_{n\ge1} \text{  weak$^*$
  converges} 
  \right\}
\end{align*}
has positive volume for both the discrete case, and the continuous
case with $f=X_1$ the time-$1$ map of a flow, where
$S_n\delta_x:= \sum_{j=0}^{n-1}\delta_{f^j x}$.  The existence of
forward Birkhoff averages for all continuous observables
$\vfi:M\to\RR$ on a full Lebesgue measure subset of points is known
for all $C^2$ globally hyperbolic (Axiom A with no cycles) systems;
see e.g. Bowen and Ruelle~\cite{Bo75,BR75}.  It is also well-known
that this is not true in general; see e.g. Ruelle and
Takens~\cite{ruelle2001,takens08}.  Moreover, as shown by Jiang, Liu
and Qian~\cite{Jiang2002}, in this setting there exists also a full
Lebesgue measure subset of Lyapunov regular points; see
e.g. Section~\ref{sec:entropy-physic-measu} for more details. In
addition, this also holds for smooth Axiom A endomorphisms and
hyperbolic SRB measures~\cite{Qian2009-rd}.

Recently, it has been shown that, on the contrary, forward Birkhoff
averages might not exist on open subsets of points for locally dense
subsets of smooth diffeomorphisms; see e.g. Kiriki, Soma, Nakano, Li
and
Vargas~\cite{KrS15,kiriki2024takens,kiriki2024pluripotency}. Moreover,
existence of Birkhoff averages on positive volume subsets of points
does not imply the existence of Lyapunov exponents on the same subset,
as shown by Kiriki, Li, Nakano, Soma~\cite{Kiriki2022}, and it may
happen that the set of Lyapunov regular points has not full Lebesgue
measure; see e.g. Ott~\cite{Ott2008}. In \cite{CCSV21}, the second
author together with Carvalho, Coelho and Varandas, proved more
general results on irregular sets with respect to C\`esaro averages,
including information on the Baire genericity of irregular sets for
averages of semigroup actions that take into account the group
structure.

Here we explore some consequences of the \emph{assumption of existence
  of forward Birkhoff averages for all continuous observables on
  positive volume subsets of points} together with \emph{existence and
  regularity of Lyapunov exponents} on the same subset of points,
motivated by the above mentioned features of Axiom A diffeomorphisms
together with the results of Tsujii~\cite{Tsujii1991}; or the
assumption of $p\omega$-\emph{regular points on positive volume
  subsets} from the results of Burguet~\cite{Burguet2019}, applied to
some robust examples of systems, both discrete and continuous.

We obtain necessary and sufficient condition on the existence of
physical measures for $3$-sectionally expanding attracting sets for
smooth flows with center-unstable \emph{mixed behavior}; and for
partially hyperbolic attractors whose center direction exhibits
non-uniformly hyperbolic \emph{mixed behavior} along the central
direction.  This means that, along certain multidimensional invariant
directions, we have both positive and negative Lyapunov exponents
Lebesgue almost everywhere.

This provides a characterization of the existence of physical measure
for the ``wild strange attractors'' for singular flows from
Shilnikov-Turaev~\cite{ST98}; and also for the $C^\infty$ smooth
instances of multidimensional non-uniformly expanding attractors from
the class of diffeomorphisms obtained by Viana~\cite{Vi97}.

We note that mostly expanding/contracting behavior along
multidimensional invariant directions is enough to ensure the
existence of physical measures, as shown by Alves, Bonatti and Viana
in~\cite{BoV00,ABV00}; but no such general results are known for mixed
behavior, among non-conservative diffeomorphisms, along
multidimensional invariant directions, apart from the case of
$C^\infty$ surface diffeomorphisms by Burguet~\cite{burguet2022srb}
(and Buzzi-Crovisier-Sarig~\cite{buzzi2022proof}), certain H\'enon
maps~\cite{BeY93}, and also certain dissipative skew-products
perturbations of the conservative Berger-Carrasco maps~\cite{BerCarr}
studied by Carrasco and Obata~\cite{obata2018,Carrobata21,obata23}.

\section{The diffeomorphism case}
\label{sec:diffeomorphism-case}

We now present preliminary definitions and results. In what follows
$M$ is a compact Riemannian manifold endowed with a Riemannian
distance $\dist$ and volume form $\m$.

\subsection{Partial hyperbolicity with central dissipativeness}
\label{sec:partial-hyperb-with}

Let $f:M\circlearrowleft$ be a diffeomorphism of $M$ having an attractor
$\Lambda$ with a \emph{trapping region} $U$ (i.e. $U$ is an open
neighborhood of $\Lambda=\Lambda_f(U):=\cap_{n\ge1}\overline{f^nU}$
with $\overline{fU}\subset U$) satisfying
\begin{enumerate}
\item \textbf{partial hyperbolicity}: the tangent bundle $T_\Lambda M$
  over $\Lambda$ splits into continuous $Df$-invariant subbundles
  $T_\Lambda M=E^s\oplus E^c\oplus E^u$ such that there are constants
  $C>0, 0<\lambda<1$ satisfying for all $x\in\Lambda$ and $n\ge1$
  \begin{enumerate}
  \item $E^s$ is uniformly contracted:
    $\|Df^n\mid_{E^s_x}\|\le C\lambda^n$;
  \item $E^u$ is uniformly expanded:
    $\|(Df^n\mid_{E^u_x})^{-1}\|\le C\lambda^n$;
  \item the splitting is dominated:
    \begin{align*}
      \max\left\{\frac{\|Df^n\mid_{ E^s_x}\|}{\|Df^n\mid_{ E^c_x}\|},
          \frac{\|Df^n\mid_{ E^c_x}\|}{\|Df^n\mid_{ E^u_x}\|}
          \right\}\le C \lambda^n;
    \end{align*}
  \end{enumerate}
\item \textbf{dissipativeness along the central direction}:
  $|\det (Df^n\mid_{ E^c_x})|\le C \lambda^n$ for all
  $x\in\Lambda, n\ge1$.
\end{enumerate}
We consider a continuous extension of each bundle of the splitting to
$U$ and assume, without loss of generality, that the extended stable
bundle is forward invariant~\cite{ArMel17}, and also that
dissipativeness and uniform expansion are preserved. We additionally
assume the following.
\begin{enumerate}[resume]
\item \textbf{Non-uniform expansion along the central
    direction}: there is $c_0>0$ such that
  \begin{align}\label{eq:wNUE}
    \limsup\nolimits_{n\to\infty}
    \log\|Df^n\mid_{ E^c_x}\|^{1/n}> c_0, \text{  for
    Lebesgue a.e.  }x\in U.
  \end{align}
\end{enumerate}

\subsection{Some subsets of regular points}
\label{sec:some-subsets-regular}

To state the main results, we recall some notions of differential
ergodic theory. As usual, we write $\log^+a=\max\{0,\log a\}$ for all
$a>0$ in what follows. For a diffeomorphism $f:M\to M$ (or for the
time-$1$ map $f=X_1$ of a flow generated by a vector field $X$), we set
the forward Lyapunov exponent of $v\in T_x M\setminus\{0\}$ as
\begin{align*}
  \chi^+(x,v):=\limsup\nolimits_{n\to\infty}\log\|Df^n(x)v\|^{1/n}.
\end{align*}
The function $\chi(x,\cdot)$ admits finitely many values only
$\chi^+_1(x)>...>\chi^+_{p(x)}(x)$ on $TM\setminus \{0\}$ and
generates a filtration
$ 0\subsetneq V_{p(x)}(x) \subsetneq \cdots \subsetneq V_{1}=T_xM$
with $V_i(x)=\{ v\in T_xM, \ \chi^+(x,v)\leq \chi^+_i(x)\}$. The
functions $p(x), \chi^+_i(x)$ are $f$-invariant; the vector spaces
$V_i(x)$, $i=1,...,p(x)$ are $Df$-invariant; and each of them depend
Borel measurably on $x$.

We write
$ \chi_1^+(x):=\limsup\nolimits_{n\to\infty}\log^+\|Df^n(x)\|^{1/n}=
\max\{\chi(x,v):v\in T_x M\setminus\{0\}\} $ for the top Lyapunov
exponent of the action of $Df:TM\to TM$ on the tangent bundle at
$x\in M$.  We also set
\begin{align*}
  \chi^+_{\wedge}(x)
  :=
  \limsup_{n\to+\infty} \frac1n \log^+\max_k\|\wedge^kDf^n(x)\|
\end{align*}
for the top $k$-sectional Lyapunov exponents of the bundle morphism
$\wedge Df$ induced by $Df$ on the exterior algebra bundle
$\wedge TM=\oplus_k\wedge^kTM$, for $k=1,\dots,d:=\dim M$; see
e.g.~\cite[Definition 1.11]{Salgado19}.

We finally set $p\omega(x)$ the family of all weak$^*$ accumulation
points of the empirical measures $\Big(\frac1n
S_n\delta_x\Big)_{n\ge1}$, which is always a nonempty set of
$f$-invariant probability measures on  compact ambient
spaces like $M$. We say that $x$ is $p\omega$-\emph{regular} if
\begin{align*}
  \sup_{\mu\in p\omega(x)}\int\chi^+_\wedge\, d\mu =\chi^+_\wedge(x).
\end{align*}
We denote by $\cP$ the subset of $p\omega$-regular points.

\subsection{Main result for diffeomorphisms}
\label{sec:main-result-diffeom}

In what follows, a \emph{hyperbolic ergodic physical measure with
  central mixed behavior} is a $f$-invariant probability measure $\mu$
such that for $\mu$-a.e. $x\in\Lambda$ and for $x$ in a positive
volume subset of $U$ we have a positive and a negative central
Lyapunov exponent: there exists $c_0>0$ so that
  \begin{align}\label{eq:mixcentral}
    \lim\nolimits_{n\to\infty}
    \log\|Df^n\mid_{ E^c_x}\|^{1/n}> c_0
           \qand
    \lim\nolimits_{n\to\infty}
    \log\|(Df^n\mid_{E^c_x})^{-1}\|^{1/n}>c_0.
  \end{align}
  Our first main result characterizes the existence of a physical
  probability measure $\mu$ with mixed central behavior.

\begin{maintheorem}\label{mthm:SRBvdiff}
  Let $f:M\to M$ be a $C^\infty$ diffeomorphism admitting a partially
  hyperbolic attracting set $\Lambda=\Lambda_f(U)$ satisfying
  conditions $(1), (2)$ and $(3)$ with a two-dimensional center
  subbundle.

  Then there exists an hyperbolic ergodic physical with mixed central
  behavior 
  if, and only if, the subset $\cP$ has positive volume in the trapping
  region $U$, i.e.  $\m(\cP\cap U)>0$.
\end{maintheorem}


\subsection{An open class of examples}
\label{sec:an-open-class}

A $C^r$ open set $\U$ of partially hyperbolic diffeomorphisms with
\emph{two-dimensio\-nal} center bundle $E^c$, in the setting of
Theorem~\ref{mthm:SRBvdiff}, is provided by Viana in \cite{Vi97} for
$r\ge3$.  Using the notation from \cite{Vi97}, we consider the
following skew-product
$\vfi_{\alpha,b}:T_3\times\RR^2\circlearrowleft,
(\Theta,X)\mapsto(\hat g(\Theta),\hat f_{\alpha,b}(\Theta,X))$ with
$\Theta=(\theta,T)\in T_3=\sS^1\times B^2$ the solid $3$-torus, where
$B^2$ is the unit $2$-disk, $X=(x,y)\in I_0^2\subset \RR^2$ and
\begin{align*}
  \hat g&:T_3\circlearrowleft, \quad \hat g(\theta,T)= (d\theta\bmod1, G(\theta,T))
\end{align*}
is a solenoid embedding; and additionally
\begin{align*}
f_{\alpha,b}&:T_3\times\RR^2\to\RR^2, \quad f_{\alpha,b}(\Theta,X)=(a(\theta)-x^2+by, -bx)
\end{align*}
where
\begin{itemize}
\item $a(\theta)=a_0+\alpha\omega(\theta)$ with
  \begin{itemize}
  \item $a_0\in(1,2)$ so that $0$ is preperiodic for the map
    $x\mapsto a_0-x^2$;
  \item $\omega:\sS^1\to\RR, \theta\mapsto \sin2\pi\theta$
    is a Morse function; and
  \end{itemize}
\item $I_0\subset (-2,2)$ is a compact interval so that
  $\vfi_{\alpha,b}(T_3\times
  I_0^2)\subset\interior(T_3\times I_0^2)$
  for every small $\alpha,b>0$.
\end{itemize}
Hence, setting $U:=T_3\times I_0^2$ we get a diffeomorphism (into its
image) $\vfi_{\alpha,b}:U\to U$ for every $b>0$. We note that
$ \det Df_{\alpha,b} = b$ and $\vfi_{\alpha,b}$ is
$C^\infty$. Moreover, clearly we have a $\vfi_{\alpha,b}$-invariant
partially hyperbolic splitting $TU=TB^2\oplus T\RR^2\oplus E^u$, where
$TB^2\oplus E^u=TT_3$ is the hyperbolic $\hat g$-invariant splitting
of the standard Smale solenoid attractor, with $E^u$ a graph
$T\sS^1\to TB^2$ over $T\sS^1$, $E^c=T\RR^2$ is the dissipative
central direction of the skew-product, and $E^s=TB^2$ the uniformly
contracting direction.

\begin{theorem}{\cite[Theorem B]{Vi97}}
\label{thm:vianadiffeo}
Let us assume $d$ to be large enough. Then, there exists an open set
of (small positive) values of $(\alpha, b)$ for which
$f:=\vfi_{\alpha,b}$ admits a partially hyperbolic attractor
$\Lambda=\Lambda_f(U)$ satisfying the conditions $(1), (2)$ and $(3)$
(with a two-dimensional center subbundle and $\liminf$ in the place of
$\limsup$ in~\eqref{eq:wNUE}).
Moreover the same holds for all $\vfi$ in a neighborhood $\U$ of
$\vfi_{\alpha,b}$ in $C^3(U)$.
\end{theorem}

\section{The flow case}
\label{sec:flow-case}

Recall that a \emph{trapping region} $U$ for a flow $X_t$ is an open
subset of the manifold $M$ which satisfies: $X_t(U)$ is contained in
$U$ for all $t>0$, and there exists $T>0$ such that
$\overline{X_t(U)} $ is contained in the interior of $U$ for all
$t>T$. We define
$\Gamma(U)=\Gamma_X(U):= \cap_{t>0}\overline {X_t(U)}$ to be the
\emph{maximal positive invariant subset in the trapping region $U$}.

A \emph{singularity} for the vector field $X$ is a point
$\sigma\in M$ such that $X(\sigma)=0$ or, equivalently,
$X_t(\sigma)=\sigma$ for all $t \in \RR$. The set formed by
singularities is the \emph{singular set of $X$} denoted
$\sing(X)$.  We say that a singularity is hyperbolic if the
eigenvalues of the derivative $DX(\sigma)$ of the vector
field at the singularity $\sigma$ have nonzero real part.

\subsection{Partial hyperbolicity}
\label{sec:partial-hyperb-secti}

A \emph{dominated splitting} over a compact invariant set $\Lambda$ of
$X$ is a continuous $DX_t$-invariant splitting
$T_{\Lambda}M = E \oplus F$ with $E_x \neq \{0\}$, $F_x \neq \{0\}$
for every $x \in \Lambda$ and such that there are positive constants
$K, \lambda$ satisfying
\begin{align*}
  \|DX_t|_{E_x}\|\cdot\|DX_{-t}|_{F_{X_t(x)}}\|<Ke^{-\la
  t}, \ \textrm{for all} \ x \in \Lambda, \ \textrm{and
  all} \,\,t> 0.
\end{align*}

A compact invariant set $\Lambda$ is said to be
\emph{partially hyperbolic} if it exhibits a dominated
splitting $T_{\Lambda}M = E \oplus F$ such that
\begin{itemize}
\item $E$ is \emph{uniformly contracted}: there are
  $C,\lambda>0$ so that $\|DX_t|_{E_x}\|\leq Ce^{-\lambda t}$ for
  $t\geq 0$ and all $x\in\Lambda$.
\end{itemize}
In this case $F$ is the \emph{central subbundle} of $\Lambda$.  Or
else, we may replace uniform contraction along $E$ by
\begin{itemize}
\item  $F$ is \emph{uniformly expanded}: there are
  $C,\lambda>0$ so that
  $\|\big(DX_t|_{F_x}\big)^{-1}\|\leq Ce^{-\lambda t}$ for $t\geq 0$
  and all $x\in\Lambda$.
\end{itemize}
We can assume, without loss of generality, that a continuous
extension of the stable direction $E$ of the splitting is still
$D\phi_t$-invariant~\cite{ArMel17} and consider also a continuous
extension of the center direction $F$ to the trapping region $U$.

\subsection{ Sectional expansion/contraction}
\label{sec:section-expans}

The following is an extension of the notion of $2$-sectional
expansion proposed in \cite[Definition 2.7]{MeMor08}.

Given $E$ a vector space, we denote by $\wedge^p E$ the exterior power
of order $p$ of $E$, defined as follows.  If $v_1,\dots, v_n$ is a
basis of $E$ then $\wedge^p E$ is generated by
$\{v_{i_1}\wedge \cdots \wedge v_{i_p}\}_{1 \leq i \leq n, i_j \neq
  i_k, j \neq k}$.  Any linear transformation $A:E\to F$ induces a
transformation $\wedge^p A:\wedge^p E\to\wedge^p F$.  Moreover,
$v_{i_1}\wedge \cdots \wedge v_{i_p}$ can be viewed as the $p$-plane
generated by $\{v_{i_1}, \cdots, v_{i_p}\}$ if
$i_j \neq i_k, j \neq k$. Reciprocally, each $p$-plane $L_x$ can be
viewed as an element $\widetilde{v}\in \wedge^p(F_x)\setminus \{0\}$
of norm one.

Given a $DX_t$-invariant subbundle $F \subset T_{\Lambda}M$ such that
$\dim F_x \geq p$ is constant for $x\in\Lambda$, we say that
\begin{itemize}
\item $F$ is  $p$-\emph{sectionally
    expanding} if 
  there exists $c,\lambda>0$ so that for all $x\in\Lambda$ and $t\ge0$
  we have 
  $\|\wedge^p \big(DX_t(x)\mid_{F_x}\big)^{-1}\|\le ce^{-\lambda t}$.
  
\item a $DX_t$-invariant subbundle
  $F \subset T_{\Lambda}M$ is \emph{mostly} $p$-\emph{sectionally
    expanding} if $\dim F_x \geq p$ is constant for $x\in\Lambda$ and
  there exists $c_0>0$ so that for Lebesgue almost every $x \in U$
  \begin{align}\label{eq:pNUE-ext}
    \limsup\nolimits_{t\to+\infty}
    \log\vert (\wedge^p DX_t \vert_{\wedge^p F_x})\vert^{1/t} > c_0.
  \end{align}
  This means that for 
every $p$-dimensional linear subspace $L_x \subset F_x$ one has
\begin{align}\label{eq:pNUE}
  \limsup\nolimits_{t\to+\infty}
  \log\vert \det (DX_t \vert_{L_x})\vert^{1/t} > c_0.
\end{align}
\end{itemize}

A \emph{mostly} $p$-\emph{sectional-hyperbolic set} is a partially
hyperbolic set whose central subbundle is mostly $p$-sectionally
expanding -- this extends the notion of sectional hyperbolicity given
by the second author in \cite[Definition 3]{Salgado19}.

This is a particular case of the so called
\emph{singular-hyperbolicity} whose definition we recall now.  A
$DX_t$-invariant subbundle $F \subset T_{\Lambda}M$ is said to be
\emph{volume expanding} if we replace~\eqref{eq:pNUE} by
\begin{align*}
  \vert \det (DX_t \vert_{F_x})\vert > C e^{\la t},
  \textrm{ for all } t>0.
\end{align*}

A \emph{singular-hyperbolic set} is a partially hyperbolic set whose
central subbundle is volume expanding~\cite[Definition 1]{MPP99}.
Clearly, in the three-dimensional case, these notions are equivalent.


\begin{remark}\label{rmk:2secexpk}
  An invariant subbundle $F$ which is $2$-sectionally expanding and
  $\dim F \geq 3$, is also $p$-sectionally expanding for each
  $p=3,\dots, \dim F$. However, the reciprocal is not true, as shown
  by the open class of examples to be discussed in what follows, since
  sectional contraction may occur in $k$-subspaces with $k<p$.
\end{remark}

To encompass  \emph{contracting behavior} along $F$, we say that
\begin{itemize}
\item an invariant subbundle $F$ of the tangent bundle $T_\Lambda M$
  is \emph{partially $k$-sectionally contracting} if there are
  $C,\lambda>0$ so that for every $x\in\Lambda$ we can find a
  $k$-dimensional subspace $L_x\subset F_x$ for which
  \begin{align}\label{eq:def-k-sec-contract}
    \vert \det (DX_t \vert_{L_x})\vert \leq C e^{-\la t},
    \textrm{ for all } t>0.
  \end{align}
\item the (extended) central-subbundle of a partially hyperbolic
  splitting is \emph{mostly $k$-sec\-tion\-al\-ly contracting} if for
  Lebesgue almost every $x\in U$, there exists a $k$-dimensional
  subspace $L_x\subset F_x$ such that
  \begin{align}
    \label{eq:pNUC}
    \limsup\nolimits_{t\to+\infty}
    \log\vert \det (DX_t \vert_{L_x})\vert^{1/t} < 0.
  \end{align}
\end{itemize}

\subsection{Main result for flows}
\label{sec:main-result-flows}

In what follows, we say that a $X_t$-invariant ergodic probability
measure $\mu$ is a \emph{hyperbolic ergodic physical with central
  mixed behavior} if, for $\mu$-a.e. $x\in\Lambda$ and for $x$ in a
positive volume subset of $U$ we have a positive and a negative
central Lyapunov exponent: there exists $c_0>0$ so that
  \begin{align*}
    \lim\nolimits_{t\to\infty}\log\|DX_t\mid_{F_x}\|^{1/t}>c_0
           \qand
    \lim\nolimits_{t\to\infty}\log\|(DX_t\mid_{F_x})^{-1}\|^{1/t}>c_0.
  \end{align*}
  We recall that $\cP$ is the subset of $p\omega$-regular points.

\begin{maintheorem}
  \label{mthm:physical3secexp}
  Let $X:M\to TM$ be a $C^\infty$ vector field admitting an attracting
  set $\Gamma_X(U)$ which is partially hyperbolic with $\dim F=3$,
  mostly $3$-sectionally expanding~\eqref{eq:pNUE} and $2$-sectionally
  contracting~\eqref{eq:pNUC}; and all singularities are hyperbolic.

  Then there exists an hyperbolic ergodic physical with mixed central
  behavior if, and only if, the subset $\cP$ has positive volume in
  the trapping region $U$, that is, we have $\m(\cP\cap U)>0$.
\end{maintheorem}

\subsection{An open class of examples}
\label{sec:an-open-class-1}

The following provides a $C^1$ robust family of vector fields in the
setting of Theorem~\ref{mthm:physical3secexp}.

\begin{figure}[htpb!]
  \centering \includegraphics[width=4.5cm]{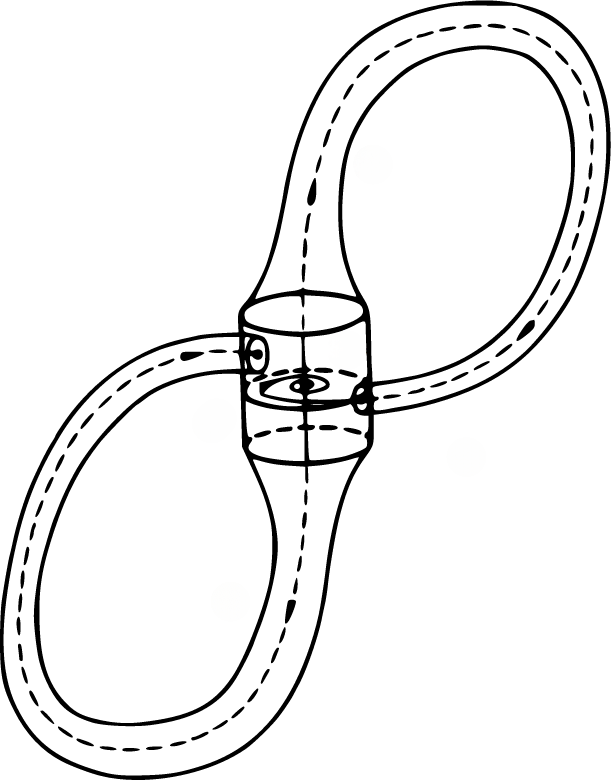}
  \caption{Sketch of the trapping region of Turaev-Shil'nikov's wild
    strange attractor around the unfolding of a saddle-focus
    equilibrium at the origin.}
  \label{fig:wild}
\end{figure}

In~\cite{ST98}, Turaev and Shilnikov construct a \emph{wild attractor}
$\Gamma$ (see Figure \ref{fig:wild}) of a vector field having a
\emph{singular hyperbolic splitting} $T_\Gamma M=E\oplus F$, where
$\dim F=3$ and $\dim E=n-3$ for any $n\ge4$; and the attractor
contains a (hyperbolic) saddle-focus singularity at the
origin. Moreover, the bundle $F$ is $2$-\emph{sectionally
  contracting}; see \cite[p 296, formula (13)]{ST98} and more
recently~\cite{Gonchenko_2021} for other concrete examples of ODE's
displaying similar features.
  
These attractors are also robust: these properties persist for a $C^1$
neighborhood $\U$ of the original vector field.  This means that the
subbundle $F$ satisfies
\begin{itemize}
\item $3$-sectional expansion, i.e.~\eqref{eq:pNUE}
  (or just the volume along $F$, since $\dim F=3$)
  holds for all $x\in U$ for $p=3$ and a fixed positive value of $t$;
  and also
\item partial $2$-sectional contraction, i.e. it
  satisfies~\eqref{eq:def-k-sec-contract}: more precisely, there exist
  $T>0$ so that for every $x\in U$ we can find a $2$-dimensional
  subspace $L_x\subset F_x$ so that $|\det (DX_T\mid_{L_x})|<1$.
\end{itemize}
Moreover, the $C^1$ robust property ensures that in $\U$ there exists
a $C^s$ dense subset of $C^r$ vector fields exhibiting the same
properties, for $2\le s \le  r \le\infty$.

\section{An abstract setting}
\label{sec:an-abstract-setting}

Now we go back to the diffeomorphism $f:M\to M$ and recall that a
point $x\in M$ is (positive) \emph{Lyapunov regular} if, for all
$1\leq k\leq d$, the \emph{sectional Lyapunov exponents} satisfy
\begin{align}\label{eq:Lyapreg}
  \chi^k(x):=
  \limsup\nolimits_{n\to+\infty}\log\|\wedge^kDf^n(x)\|^{1/n}
  =\sum\nolimits_{i=1}^k\chi^+_i(x).
\end{align}
By Oseledets' Theorem~\cite{Os68}, the set $\cL$ of Lyapunov regular
points has full measure for any $f$-invariant probability measure.

Motivated by Tsujii~\cite{Tsujii1991}, we say that $x\in M$ is
$T$-\emph{regular} for $f$ if $x\in\BB$ (so that $p\omega(x)=\{\mu_x\}$)
and $x$ is $p\omega$-regular, and also Lyapunov
regular~\eqref{eq:Lyapreg} for $\mu_x$, that is, we have
$\chi^+_i(x)=\chi^+_i(y)$ for $\mu_x$-a.e. $y$.

We denote by $\cR=\BB\cap\cL\cap\cP$ the subset of $T$-regular points
in what follows.

\subsection{Unimodal sectional Lyapunov spectrum and existence of
  physical measures}
\label{sec:unimod-section-lyapu-1}

We say that the \emph{sectional Lyapunov spectrum} of a point
$x\in M$, given by $\chi^i(x), i=1,\dots,p(x)$, is \emph{unimodal} if
there exists $1\le k(x)<p(x)$ so that
\begin{align}
  \label{eq:unimodal}
  \chi^{p(x)}(x)<\cdots<\chi^{k(x)}(x)
  \qand
  \chi^{k(x)}(x)>\dots>\chi^1(x).
\end{align}
If in addition $x\in\cR$, then the Lyapunov regularity ensures that
the $\mu_x$ is hyperbolic, since 
\begin{align*}
  \chi^+_1(x)=\chi^1(x) \qand
  \chi^+_i(x)=\chi^{i}(x)-\chi^{i-1}(x), \quad i=2,\dots, d=\dim M;
\end{align*}
and thus
$\chi^+_d(x)\le\cdots\le\chi^+_{k(x)+1}(x)<0<\chi^+_{k(x)}\le\cdots\le\chi^+_1(x)$.

When $f=X_1$ is the time-$1$ map of a flow, then we allow for the zero
Lyapunov exponent along the flow direction: we say that $x$ has
\emph{unimodal sectional Lyapunov spectrum} for a flow if there
exists $1< k(x)<p(x)$ so that
\begin{align}
  \label{eq:unimodalflow}
  \chi^{p(x)}(x)<\cdots<\chi^{k(x)}(x)
  \,\&\,
  \chi^{k(x)}(x)=\chi^{k(x)-1}(x)>\dots>\chi^1(x).
\end{align}


The condition~\eqref{eq:unimodal} or~\eqref{eq:unimodalflow} is
enough to obtain a physical/SRB measure for $C^\infty$ maps or flows
\emph{assuming simultaneously that the collection of such points has
  positive volume inside the family $\cR$ of regular points}. The
reciprocal also holds true as a consequence of the work of Jiang et
al~\cite{Jiang2002}.

\begin{maintheorem}
  \label{mthm:C2physical}
  Let $f:M\to M$ be a $C^\infty$ diffeomorphism (e.g. the time-$1$ map
  of a flow $X_t$ induced by a $C^\infty$ vector field).  Then there
  exists an invariant \emph{hyperbolic ergodic physical} probability
  measure if, and only if, the subset $S$ of points with unimodal
  spectrum has positive volume among the regular points $\cR$, i.e.
  $\m(S\cap\cR)>0$.
\end{maintheorem}

\subsection{Open classes of examples}
\label{sec:open-classes-satisfy}

Among vector fields, a class that satisfies the assumptions of
Theorem~\ref{mthm:C2physical} are the sectional-hyperbolic attracting
sets for $C^2$ vector fields, which admit physical measures; see
e.g~\cite{araujo_2021} and references therein. This includes, in
particular, hyperbolic attractors for $C^2$ vector fields~\cite{BR75}
and all Lorenz-like vector fields~\cite{AraPac2010s}.

Among diffeomorphisms, a class that satisfies the assumptions of
Theorem~\ref{mthm:C2physical} is the family of non-uniformly
hyperbolic attracting sets for $C^2$ diffeomorphisms with
H\"older-continuous splitting which, as shown in~\cite{CMY2017}, admit
physical measures. This class includes non-uniformly partially
hyperbolic attractors, as in~\cite{ABV00}, as well as hyperbolic
attractors and Anosov diffeomorphisms.

For these classes we may specialize the statement of
Theorem~\ref{mthm:C2physical} as follows.

\begin{maincorollary}
\label{mcor:NUHav}
Let $f:M\to M$ be a $C^\infty$ diffeomorphism with an attracting
subset $\Lambda$, admitting a trapping region $U$ with a continuous
invariant splitting $TU=E^{cs}\oplus E^{cu}$, so that on a positive
Lebesgue measure subset $\Omega\subset U$ we have $c_0>0$ so that for each
$x\in\Omega$
  \begin{enumerate}
  \item
    $ \liminf_{n\nearrow\infty}\frac1n
    \sum_{j=0}^{n-1}\log\|Df\mid_{E^{cs}_{f^jx}}\| < -c_0$, i.e.,
    $E^{cs}$ is non-uniformy contracting (on average); and
  \item
    $\liminf_{n\nearrow\infty}\frac1n \sum\nolimits_{j=0}^{n-1}
    \log\big\|\big(Df\mid_{E^{cu}_{f^jx}}\big)^{-1}\big\| < -c_0$,
    i.e., $E^{cu}$ is non-uniformy expanding (on average).
  \end{enumerate}
  If $\m(\cR\cap\Omega)>0$, then there exists an ergodic hyperbolic
  invariant physical probability measure.
\end{maincorollary}

\section{Comments and conjectures}
\label{sec:comments-conjectures}

The statements Theorems~\ref{mthm:SRBvdiff},
~\ref{mthm:physical3secexp} and~\ref{mthm:C2physical} are corollaries
of the results of Burguet~\cite{Burguet2019}, together with the
results from Jiang et al. from~\cite{Jiang2002} on the
extension of Lyapunov regularity along stable manifolds, and the
definition of regular point is inspired by the main results from
Tsujii~\cite{Tsujii1991}. The proofs are detailed in
Section~\ref{sec:entropy-physic-measu}.

\begin{remark}[Regularity, no hyperbolicity and no physical measure]
  \label{rmk:hypness}
  We can have a full Lebesgue measure subset of $p\omega$-regular and
  Lyapunov regular points without physical measure: consider the
  following skew-product on the $3$-torus
  \begin{align*}
    T:\TT^3\to\TT^3, \quad (x,y,z)\mapsto (x+y,y, h(z))
  \end{align*}
  where $h:\sS^1\circlearrowleft$ is the ``North-South map''; see
  Figure~\ref{fig:NSmap}.
  \begin{figure}[htpb]
    \centering
    \includegraphics[width=3cm]{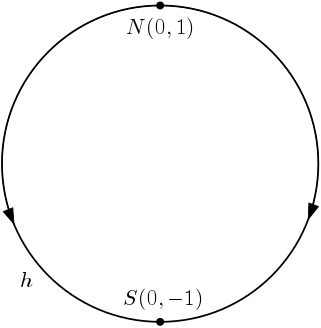}
    \caption{The North-South map $h:\sS^1\to\sS^1$ in the circle.}
    \label{fig:NSmap}
  \end{figure}
  The restriction $T_0:=T\mid_{\TT^2\times\{S\}}$ to the invariant
  $2$-torus $\TT^2\times\{S\}$ is conservative and each circle
  $\sS^1\times\{y\}\times\{S\}$ is $T_0$-invariant for each given
  $y\in\sS^1$. Moreover, the dynamics on these circles is a rotation:
  $T_0(x,y,S)=(x+y,y,S)$, which preserves the Lebesgue length measure
  $\lambda$ on each cicle. Hence, $\lambda$-a.e. points are Lyapunov
  and $p\omega$-regular. Since this holds for all $y$, then we have
  that a.e. points $(x,y,S)$ are also regular with respect to the
  Lebesgue area measure on $\TT^2\times\{S\}$.

  In addition, all points $(x,y,z)\in \TT^3$ with
  $z\in\sS^1\setminus\{N\}$ belong to the stable manifold of
  $(x,y,S)$. Thus, from Theorem~\ref{thm:ErgLyapExp}, we have that
  a.e. point of $\TT^3$ are Lyapunov and $p\omega$-regular.

However, there is no ergodic physical measure for $T$, since all
ergodic measures are the Lebesgue length measures along the circles
$\sS^1\times\{y\}\times\{S\}$ for $\lambda$-a.e. $y\in\sS^1$
(corresponding to irrational rotations), whose ergodic basins are
given by $\sS^1\times\{y\}\times(\sS^1\setminus\{N\})$, with zero
volume on $\TT^3$.

Finally, we note that the system is not hyperbolic: the Lyapunov
exponents of $T_0$ are all null.
\end{remark}

Corollary~\ref{mcor:NUHav} yields a non-trivial result since, under
the conditions of items (1) and (2), it is known how to construct a
$cu$-Gibbs state (or a hyperbolic $SRB$ measure)
following~\cite{ABV00}. However, in general, these measures are not
necessarily hyperbolic, as the following example shows.

\begin{remark}[Existence of non-hyperbolic ergodic $cu$-Gibbs state]
  \label{rmk:nonsuperfl} 
  Indeed, we consider a pair of diffeomorphisms
  $f,g:\cS^1\times\DD\circlearrowleft$, where $f$ is the uniformly
  hyperbolic Smale solenoid map, see
  e.g.~\cite[Sec. 7.7]{robinson1999}; and $g$ its ``intermittent''
  modification~\cite[Sec. 4.6]{Alves2020b} from~\cite{AlPi08}. In both
  cases we have attractors (i.e. transitive attracting sets)
  $\Lambda_f, \Lambda_g$ with partially hyperbolic splitting
  $E^s\oplus E^{cu}$ and ergodic (in fact, mixing) hyperbolic
  $cu$-Gibbs states $\mu_f,\mu_g$ which are the unique physical
  measures. However, $f$ is uniformly hyperbolic, while $g$ admits a
  fixed point $p\in\Lambda_g$ so that $Dg_p\mid_{E^{cu}}$ is an
  isometry.

  Hence, if $\ell>1$ is such that expansion/contraction rates of
  $f^\ell$ are stronger than the ones of $g$, then $F:=f^\ell\times g$
  is partially hyperbolic and $\nu=\mu_f\times\mu_g$ is an ergodic
  $F$-invariant measure which is the unique physical measure and a
  $cu$-Gibbs state. Thus, $F$ is nonuniformly hyperbolic on a full
  volume measure subset of $(\cS^1\times\DD)^2$ and the attractor
  $\Lambda:=\Lambda_f\times\Lambda_g$ admits the partially hyperbolic
  splitting
  $T_\Lambda M=\big(E^s_f\oplus E^s_g\big)\oplus \big(E^{cu}_g\oplus
  E^{cu}_f\big)$. However, the ergodic $cu$-Gibbs state
  $\nu=\mu_f\times\delta_p$ is non-hyperbolic, with a zero Lyapunov
  exponent along the direction $E^{cu}_g$. This shows that \emph{even
    with a full volume of non-uniformly hyperbolic points and unique
    physical/SRB measure, there can be ergodic $cu$-Gibbs states which
    are not hyperbolic}.
\end{remark}

We note that the central bundle dynamics of the skew-products
$\vfi_{\alpha,b}$ in the statement of Theorem~\ref{thm:vianadiffeo} is
a perturbation of the H\'enon family~\cite{He76}. For this polynomial
family of dissipative maps, a positive Lebesgue measure of parameters
was obtained by Benedicks and Carleson~\cite{BC91} for which the
system is non-uniform expanding, that is, condition $(3)$ is
satisfied. The recent results for $C^\infty$ smooth surface dynamics
of Burguet~\cite{burguet2022srb} would automatically ensure the
existence and uniqueness of a physical measure for these system, which
was obtained with a completely different strategy, specifically adapted
to these systems, by Benedicks and Young~\cite{BeY93}.

The conditional existence results of physical measures stated in
Theorems~\ref{mthm:SRBvdiff} and~\ref{mthm:physical3secexp}, together
with the positive results on ``Takens last problem''~\cite{KrS15},
naturally raise the following questions.

\begin{conjecture}
  \label{conj:existBirkhoff}
  For a dense/generic subset of dynamical systems, there exists a
  positive volume subset in the ambient space where at least one of
  the following properties hold:
  \begin{enumerate}
  \item the points are $p\omega$-regular;
  \item forward Birkhoff averages exist on all points;
  \item forward Lyapunov exponents exist and are Lyapunov regular.
  \end{enumerate}
\end{conjecture}

This would automatically provide for the existence of physical
measures among the smooth systems in the open classes of non-uniformly
hyperbolic flows and diffeomorphims described in the previous
subsection, and provide for a criterion to obtain physical/SRB
measures using the sectional Lyapunov spectrum.

Since Hertz, Hertz, Tahzibi and Ures~\cite{HHTU11} have shown the
uniqueness of SRB measures for transitive diffeomorphisms on surfaces,
and Obata~\cite{obata23} obtained uniqueness for non-conservative
perturbations of certain skew-products with two-dimensional center
bundle with mixed behavior, the fact that the non-hyperbolic center
direction of the family $\vfi_{\alpha,b}$ is two-dimensional suggests
the following.

\begin{conjecture}
  \label{conj:uniqSRB}
  There exists at most one SRB/physical measure for each of the maps
  $\vfi_{\alpha,b}$ in the setting of Theorem~\ref{mthm:SRBvdiff}.
\end{conjecture}

Assuming a positive answer to the previous conjectures, it is natural
to consider the dependence of the physical measure on the dynamics of
the open families $\U$ of maps considered in
Subsections~\ref{sec:an-open-class}
and~\ref{sec:an-open-class-1}. Ideally we would like to obtain
\emph{statistical stability}, that is, if $\mu_n$ is the physical
measure for the dynamics of $f_n$, and $f_n$ converges to $f\in\U$ in
some $C^r$ topology ($2\le r\le\infty$), then $\mu_n$ weakly converges
to a $f$-invariant probability measure $\mu$, which is also a physical
measure for $f$.

\begin{conjecture}
  \label{conj:statstab}
  The physical measures of the $C^\infty$ dynamics in $\U$ are
  statistically stable.
\end{conjecture}

This would naturally lead to obtain physical measures for most
dynamics among the open classes of systems considered here. A natural
extension of this work is to study statistical properties (e.g. mixing
and speed of mixing) of these physical measures.

\subsection*{Acknowledgments}

V.A. thanks the Mathematics and Statistics Institute of the Federal
University of Bahia (Brazil) for its support of basic research and
CNPq (Brazil) for partial financial support. L.S. thanks the
Mathematics Institute of Universidade Federal do Rio de Janeiro
(Brazil) for its encouraging of mathematical research, and the
Mathematics and Statistics Institute of the Federal University of
Bahia (Brazil) for its hospitality, as well as FAPERJ-Funda\c c\~ao
Carlos Chagas Filho de Amparo \`a Pesquisa do Estado do Rio de Janeiro
Projects APQ1-E-26/211.690/2021 SEI-260003/015270/2021 and
JCNE-E-26/200.271/2023 SEI-260003/000640/2023, by Coordena\c c\~ao de
Aperfei\c coamento de Pessoal de N\'ivel Superior CAPES — Finance Code
001 and PROEXT-PG project Dynamic Women - Din\^amicas, CNPq-Brazil
(grant Projeto Universal 404943/2023-3), for partial financial
support.

\section{Physical measures for smooth systems}
\label{sec:entropy-physic-measu}

We present first the main result from Burguet~\cite{Burguet2019},
involving the subset $\cP$, needed to obtain the main statements.

\begin{theorem}{\cite[Main Theorem]{Burguet2019}}
  \label{thm:burguet}
  Let there be given a $C^\infty$ map $f:M\to M$. Then for Lebesgue
  almost every point $x$ there exists $\mu_x\in p\omega(x)$ such that
  $h_{\mu_x}(f)\ge \chi_\wedge^+(x)$.
\end{theorem}

We note that, in the setting of Theorem~\ref{thm:burguet}, we obtain
for $\m$-a.e. $x\in\cP$
\begin{align*}
  \chi^+_\wedge(x)
  \le
  h_{\mu_x}(f)
  \le
  \int \chi^+_\wedge \,d\mu_x
  =
  \chi^+(\mu)
  \le
  \sup_{\mu\in p\omega(x)}\int \chi^+_\wedge \,d\mu
  =\chi^+_\wedge(x),
\end{align*}
where the last inequality is given by $p\omega$-regularity of $\mu_x$,
and we used the \emph{Ruelle Inequality}~\cite{Ru78} for the entropy
$ h_{\nu}(f)\le\int \chi_\wedge^+ \, d\nu $ for any $f$-invariant
probability measure $\nu$.  Therefore, $\mu_x$ satisfies the
\emph{Pesin Entropy Formula} $h_{\mu_x}(f)=\int\chi^+_\wedge\,d\mu_x$.

In the setting of the Main Theorems, we have at least a positive
Lyapunov exponent, thus we can find an ergodic component of $\mu_x$
still satisfying the Entropy Formula with positive entropy -- we keep
the notation $\mu_x$ in what follows.

We recall that the stable (unstable) manifold of a point $x\in M$ is
the subset of points whose future (past) orbits approach the orbit of
$x$ exponentially fast:
\begin{align*}
  W^s_x
  &:=
    \{ y\in M:
  \limsup\nolimits_{n\to+\infty}
  (1/n) \log \dist\big( f^n(y), f^n(x)\big)<0 \}; \qand
  \\
  W^u_x
  &:=
    \{ y\in M:
  \limsup\nolimits_{n\to-\infty}
  (1/n)\log\dist\big( f^n(y), f^n(x)\big)<0\}.
\end{align*}
It it well-known that, under some hyperbolicity conditions, such
subsets are immersed differentiable submanifolds of $M$; see
e.g.~\cite{BarPes2007}. In what follows, given a subset $V\subset M$,
we write $W^s(V)=\cup\{ W^s_x: x\in V\}$.

From the results of Ledrappier and Young~\cite{LY85},
the measures $\mu_x$ satisfying the Pesin Entropy Formula with
positive entropy are Sinai-Ruelle-Bowen (SRB) measures, which are
dynamically characterized as \emph{$\mu_x$-almost every point has a
  positive Lyapunov exponent and the disintegration of $\mu$ along the
  unstable manifolds is absolutely continuous with respect to the
  volume on the unstable manifolds inherited from the Riemannian
  structure on the ambient manifold}.

From Pesin's Nonuniformly Hyperbolic Theory, it is known that any
ergodic hyperbolic SRB measure is physical~\cite{Pe77,BarPes2007}.
Consequently if, in addition, $\mu_x$ has no zero Lyapunov exponents
(or only one zero Lyapunov exponent along the direction of the flow),
then $\mu_x$ becomes a physical measure.

Next we present a result from Jiang et al~\cite{Jiang2002} enabling us
to show that existence of physical/SRB measures implies Lyapunov and
$p\omega$-regularity for a subset of point with positive volume; see
also~\cite[Chapter 8]{JiangQiPing03}.

\begin{theorem}{\cite[Theorem 5.3]{Jiang2002}}
\label{thm:ErgLyapExp}
Let \( f:M\to M \) be a \( C^2 \) diffeomorphism and \( \mu \) an
\( f \)-ergodic probability measure on \( M \) satisfying:
  \begin{enumerate}[(i)]
  \item $\mu$ is hyperbolic, that is, \( f \) has almost no zero
    Lyapunov exponents with respect to \( \mu \);
  \item \( \mu \) has \emph{absolutely continuous conditional
      measures on unstable manifolds} or, equivalently, \emph{Pesin’s
      entropy formula} $h_\mu(f)=\int\chi^+_\wedge\,d\mu$ holds true
    for the system \( (M, f, \mu) \).
  \end{enumerate}
  Let \( \cL=\cL(\mu) \) be the set of Lyapunov regular points
  which are contained in \( B(\mu) \). Then \( \m(W^s(\cL)) > 0 \) and
  \( W^s(\cL) 
  \subset \cL\cap B(\mu) \).
\end{theorem}

That is, each point \( y \in W^s(\cL) \) is Lyapunov regular:
$\chi^+_i(y)=\chi^+_i(x)$ and $\chi^k(y)=\chi^k(x)$, for $y\in W^s_x$
and $x\in\cL$. 
Moreover, for $\mu$-a.e. $x$ we also get
$\mu= \lim_{n \to +\infty} \frac{1}{n} \sum_{k=0}^{n-1} \delta_{f^k
  y}$.

\begin{remark}[Lyapunov and $p\omega$-regularity]
  \label{rmk:Lyapwreg}
  This ensures that $y\in W^s_x$ is also $p\omega$-regular for
  $x\in\cL\cap B(\mu)$. Indeed, since $\mu$ is $f$-invariant and
  ergodic, then $p\omega(y)=\{\mu\}$ and
  $h_\mu(f)=\int \chi^+_\wedge\,d\mu=\chi^+(x)$ for $\mu$-a.e. $x$,
  which is contained in $B(\mu)\cap\cL$, and also
  $\chi^+_\wedge(y)=\chi^+_\wedge(x)$ from
  Theorem~\ref{thm:ErgLyapExp}. Thus $h_\mu(f)=\chi^+(y)$. Comparing
  with~\cite[Remark 5]{Burguet2019} on ``Bowen's figure 8 example''
  and Remark~\ref{rmk:hypness} we see that condition $(ii)$ in the
  statement of Theorem~\ref{thm:ErgLyapExp} is
  key to obtain positive volume of $p\omega$-regular points.
\end{remark}

Hence, we only need to prove the necessary condition in the statements
of the main theorems, since the existence of a physical hyperbolic
ergodic measure $\mu$, which also satisfies the entropy formula,
automatically ensures $p\omega$-regularity and $T$-regularity for a
positive Lebesgue measure subset inside the ergodic basin of $\mu$.

\subsection{Proof of the diffeomorphism case}
\label{sec:proof-diffeom-case}

A consequence of the assumptions of Theorem~\ref{mthm:SRBvdiff} and
Theorem~\ref{thm:burguet}, is that for a $C^\infty$ map $f$ and
$x\in U\cap\cP$, we have that $\mu_x$ is an SRB measure. To conclude
the statement of Theorem~\ref{mthm:SRBvdiff} it is enough to show that
$\mu_x$ is hyperbolic with mixed behavior along the central direction.

Let then $\mu_x$ be an SRB measure with $x\in U$
satisfying~\eqref{eq:wNUE}. Since
$T_\Lambda M=E^s\oplus E^c\oplus E^u$ and $E^s$ is uniformly
contracting, $E^u$ is uniformly expanding, and $E^c$ is
two-dimensional and area contracting, then for
$1\le i\le d_u:=\dim E^u$ and $n\ge1$ we get
\begin{align}\label{eq:Eu}
  \|\wedge^i Df^n\|
  \ge
  \|\wedge^i (Df^n\mid_{E^u_x})\|
  &\ge
  \|\wedge^i (Df^n\mid_{E^u_x})^{-1}\|^{-1}
  \ge
  \|(Df^n\mid_{E^u_x})^{-1}\|^{-i}
  \ge
  C^{-i}\lambda^{-ni}
\end{align}
and so
$\chi^{d_u}(x)\ge\chi^{d_u-1}(x)\ge\cdots\ge\chi^1(x)>\log\lambda^{-1}>0$.
From domination we get, for each $u\in\wedge^{d_u}E^u_x$ and $v\in E^c_x$
\begin{align*}
  \|\wedge^{d_u+1}Df_x^n\|
  \ge
  \|(\wedge^{d_u+1}Df_x^n)(v\wedge u)\|
  \approx
  \|Df^n_xv\| \cdot\|\wedge^{d_u}Df^n_x u\|,
\end{align*}
where $\approx$ means equality up to a bounded multiplicative
constant. From~\eqref{eq:wNUE} we can find $v\in E^c_x$ so that
$\chi^+(x,v)>c_0>0$ and from~\eqref{eq:Eu} we can find
$u\in \wedge^{d_u}E^u_x$ which realizes the rate of volume growth
along $E^u$, so we obtain
$\chi^{d_u+1}(x)=\chi^+_{d_u+1}(x)+\chi^{d_u}(x)>\chi^{d_u}(x)+c_0>\chi^{d_u}(x)$. Moreover,
\begin{align*}
  \|\wedge^{d_u+2}Df_x^n\|
  &\approx
  \|\wedge^2(Df^n\mid_{E^c_x})\|\cdot\|\wedge^{d_u}(Df^n\mid_{E^u_x})\|
  \\
  &=
    |\det(Df^n\mid_{E^c_x})|\cdot\|\wedge^{d_u}(Df^n\mid_{E^u_x})\|
    \le C\lambda^n \|\wedge^{d_u}(Df^n\mid_{E^u_x})\|,
\end{align*}
and therefore we conclude that
$
\chi^{d_u+2}(x)<\chi^{d_u}(x)+\log\lambda<\chi^{d_u}+c_0<\chi^{d_u+1}(x).
$

In addition, since $E^s$ is uniformly contracting, then we also get
$\chi^{d_u+2+i+1}(x)<\chi^{d_u+2+i}(x)$ for all $i=0,\dots, \dim E^s$.
This ensures that $\chi^{d_u+1}(x)=\chi^+_\wedge(x)=h_{\mu_x}(f)$.
Thus, $\mu_x$ is hyperbolic with a Lyapunov spectrum satisfying
\begin{align*}
 \chi^+_{d_u+2}(x)=\chi^{d_u+2}(x) - \chi^{d_u+1}(x)<\log\lambda-c_0<0<c_0<\chi^+_{d_u+1}(x)
\end{align*}
Hence, the minimal exponent along the subbundle $E^c_x$ is smaller
than $\log\lambda-c_0$.  This concludes the proof of
Theorem~\ref{mthm:SRBvdiff}.

\subsection{Proof of the flow case}
\label{sec:proof-flow-case}

Similarly, to conclude the statement of
Theorem~\ref{mthm:physical3secexp} it is enough to show that $\mu_x$
is hyperbolic with mixed behavior along the central direction for
$\m$-a.e.  $x\in U\cap\cP$ and a $C^\infty$ vector field $X$. Then, from
$3$-sectional expansion, we obtain
\begin{align*}
  \chi_3^+(x)+\chi_2^+(x)+\chi^+_1(x)\ge\chi^3(x)>0,
\end{align*}
which ensures that $\chi^+_1(x)>0$.  Since the hyperbolic
singularities within $U$ are hyperbolic and are fixed points of the
dynamics, then they belong to $\Gamma$ and the sectional expansion
assumption ensures that each singularity is of saddle-type. Hence, the
union of the stable manifold of the singularities has zero
volume\footnote{Since the stable manifold of a hyperbolic critical
  element in an immersed submanifold \cite{robinson1999}, then it has
  zero volume as a subset of the ambient manifold}, and so $\m$-a.e.
$x\in U$ has a zero Lyapunov exponent along the flow
direction. Moreover, the flow direction is contained in $F$ due to the
following.
\begin{lemma}{\cite[Lemma 3.2]{arsal2012a}}
  \label{le:flow-center}
  Let $\Lambda$ be a compact invariant set for the vector field $X$.
  \begin{enumerate}
  \item Given a continuous splitting $T_\Lambda M =
    E\oplus F$ such that $E$ is uniformly contracted,
    then $X(x)\in F_x$ for all $x\in \Lambda$.
  \item Assuming that $\Lambda$ is non-trivial and has a continuous
    \emph{and dominated} splitting $T_\Lambda M = E\oplus F$ such that
    $X(x)\in F_x$ for all $x\in\Lambda$, then $E$ is a uniformly
    contracted subbundle.
  \end{enumerate}
\end{lemma}
In addition, from $2$-sectional contraction, we get that there exists
$v\in F_x\setminus\{0\}$ so that $\chi^+(x,v)<0$, and thus we conclude
that $\chi^+_3(x)<0=\chi^+_2(x)$.


This shows that $\chi^+_\wedge(x)=\chi^1(x)=\chi^2(x)$, that
$\mu=\mu_x$ is hyperbolic (since $E$ is uniformly contracted) and,
consequently, a \emph{physical measure} with mixed central
behavior for the time-$1$ map of the flow.

The probability measure $\nu:=\int_0^1(X_t)_*\mu\,dt$ is
$X_t$-invariant for every $t>0$ and also ergodic for the flow.
Moreover, if $\vfi:M\to\RR$ is continuous and $x\in B(\mu)$, then
$\psi=\int_0^1\vfi\circ X_s\,ds$ is also continuous, and since
$X_t$ and $f$ commute
  \begin{align*}
    \int\vfi\,d\nu
    &=
    \int\psi\,d\mu
    =
    \lim_{n\to+\infty}\frac1n\sum\nolimits_{j=0}^{n-1}\psi(f^j(x))
    \\
    &=
    \lim_{n\to+\infty}\frac1n\int_0^n\vfi(X_s(x))\,ds
    =
    \lim_{T\nearrow\infty, T\in\RR}\frac1T\int_0^T\vfi(X_s(x))\,ds,
  \end{align*}
  where the last equality follows from boundedness of $\vfi$. This
  shows that $B(\mu)\subset B(\nu)$ and so $\nu$ is a
  physical measure for the flow.

  This completes the proof of Theorem~\ref{mthm:physical3secexp}
  since, for $\mu$-a.e. $y\in B(\mu)\cap U\subset B(\nu)\cap U$
\begin{align*}
  \lim_{t\to\infty}\log\|DX_t\mid_{F_y}\|^{1/t}
  &=
  \limsup_{n\to+\infty}\log\|DX_1^n\mid_{F_y}\|^{1/n}
  = \chi^+_1(y)>0 \quad\text{and also}
\\
  \lim_{t\to\infty}\log\|(DX_t\mid F_y)^{-1}\|^{1/t}
  &=
  \limsup_{t\to\infty}\log\|(DX_1^n\mid_{F_y})^{-1}\|^{1/n}
  =
  -\chi^+_3(y)>0.
\end{align*}

\begin{remark}\label{rmk:noregular}
  Up to this point in this section we have not assumed that $x$ is
  Lyapunov regular.
\end{remark}

\subsection{Proof of the abstract result}
\label{sec:proof-abstract-resul}

To prove Theorem~\ref{mthm:C2physical}, we note that, by assumption, we
can find $x\in\cR$ which is simultaneously $p\omega$-regular, Lyapunov
regular and has unimodal sectional Lyapunov spectrum, so that the
measure $\mu_x\in p\omega(x)$, given by Theorem~\ref{thm:burguet},
satisfies Pesin's Entropy Formula
$h_{\mu_x}(f)=\int\chi^+_{\wedge}\,d\mu_x$.

From Ergodic Decomposition~\cite[Chapter II, Section 6]{Man87} and
Jacobs' Theorem~\cite[Theorem 9.6.2]{ViOl16} we can assume without
loss of generality that $\mu_x$ is ergodic (otherwise, there exists an
ergodic component which also satisfies Pesin's Entropy Formula) and is
hyperbolic (since the Lyapunov exponents $\chi^+_i(x)$ coincide with
the exponents $\chi^+_i(y)$ for $\mu_x$-a.e. $y$).
Hence, $\mu_x$ is a hyperbolic ergodic physical measure as claimed.
This completes the proof of Theorem~\ref{mthm:C2physical}.

For Corollary~\ref{mcor:NUHav}, we just have to note that for
$x\in \Omega\cap\cR$ the following limits exist and are bounded as
follows
\begin{align*}
  \lim_{n\nearrow\infty}\frac1n\log\|Df^n\mid_{E^{cs}_x}\|
  &\le
    \liminf_{n\nearrow\infty}\frac1n
    \sum\nolimits_{j=0}^{n-1}\log\|Df\mid_{E^{cs}_{f^jx}}\| <-c_0;
    \qand
  \\
  \lim_{n\nearrow\infty}\frac1n \log\big\|(Df^n\mid_{E^{cu}_x})^{-1}
  \big\|
  &\le
    \liminf_{n\nearrow\infty}\frac1n \sum\nolimits_{j=0}^{n-1}
    \log\|(Df\mid_{E^{cu}_{f^jx}})^{-1}\| < -c_0.
\end{align*}
This ensures that points in $\Omega\cap\cR$ have unimodal
spectrum. Thus, the system satisfies the conditions of
Theorem~\ref{mthm:C2physical} on a positive Lebesgue measure subset,
and the conclusion follows.




\def\cprime{$'$}

 \bibliographystyle{abbrv}

\end{document}